\documentclass[12pt,onecolumn]{IEEEtran}
\usepackage{setspace}
\usepackage{amssymb}
\usepackage{cite}
\usepackage{graphicx}
\usepackage{tikz}
\usepackage{pgfplots}
\usepackage{enumitem}
\usepackage{color,soul}
\usepackage[cmex10]{amsmath}
\usepackage{fancyhdr}
\ifCLASSOPTIONcompsoc
  \usepackage[caption=false,font=normalsize,labelfont=sf,textfont=sf]{subfig}
\else
  \usepackage[caption=false,font=footnotesize]{subfig}
\fi
\usepackage{fixltx2e}
\renewcommand\thesection{\arabic{section}}
\renewcommand\thesubsection{\thesection.\arabic{subsection}}

\begin{document}

\begin{flushleft}
\textbf{\textsf{\large{Augmenting Numerical Stability of the Galerkin Finite Element Formulation for Electromagnetic Flowmeter Analysis}}}

Sethupathy S.$^{1*}$, Udaya Kumar$^{2}$ \\
$^{1, 2}$ Department of  Electrical Engineering, Indian Institute of Science, Bangalore 560012, India
$^*$sethupathys87@gmail.com
\end{flushleft}

\vfill
This paper is a preprint of a paper accepted by \textit{IET Science, Measurement \& Technology} and is subject to Institution of Engineering and Technology Copyright. When the final version is published, the copy of record will be available at IET Digital Library.

\newpage

\textbf{ Abstract: Due to the complexities in handling liquid metals, theoretical evaluation of the sensitivity of magnetic flow meters forms an attractive and preferred choice. The classical Galerkin finite element formulation is generally opted for the required evaluation.  However, it is known to lead to numerical oscillations at higher flow rates. To overcome this, modified methods like upwind/Petrov-Galerkin schemes are generally suggested in allied areas like fluid dynamics. However, it requires the evaluation of stabilization parameter and this parameter is not readily available for elements of order beyond quadratic. After a careful analysis of the numerical instability through a reduced one dimensional problem, an elegant and stable approach is devised. In this scheme, the input magnetic field is restated in terms of the associated vector potential and the classical {Galerkin Finite Element Method (GFEM)} is employed without any modification. The analytical solution of the associated difference equation is employed to show: (i) the stability of the proposed approach at higher flow rates and (ii) quantification of the small oscillations what it introduces at intermediate flow rates.  It is then applied to the original  flowmeter problem and the stability of the numerical solution is clearly demonstrated.}

\IEEEpeerreviewmaketitle

\section{\bf{Introduction}}
Electromagnetic flowmeter is a non-invasive instrument which is widely used in fast-breeder reactors for the measurement of flow rate of liquid metals. As an accurate measurement of flow rate 
is essential for the safe operation and control of the reactor, the performance of flowmeter needs to be reliably ascertained. Due to the practical difficulties in handling liquid metals, experimental determination of flowmeter sensitivity is an involved job. 
Hence, accurate theoretical or numerical evaluation of the sensitivity formed an attractive alternative.

\begin{figure}[htp]
		\centering
		\includegraphics[scale=0.4]{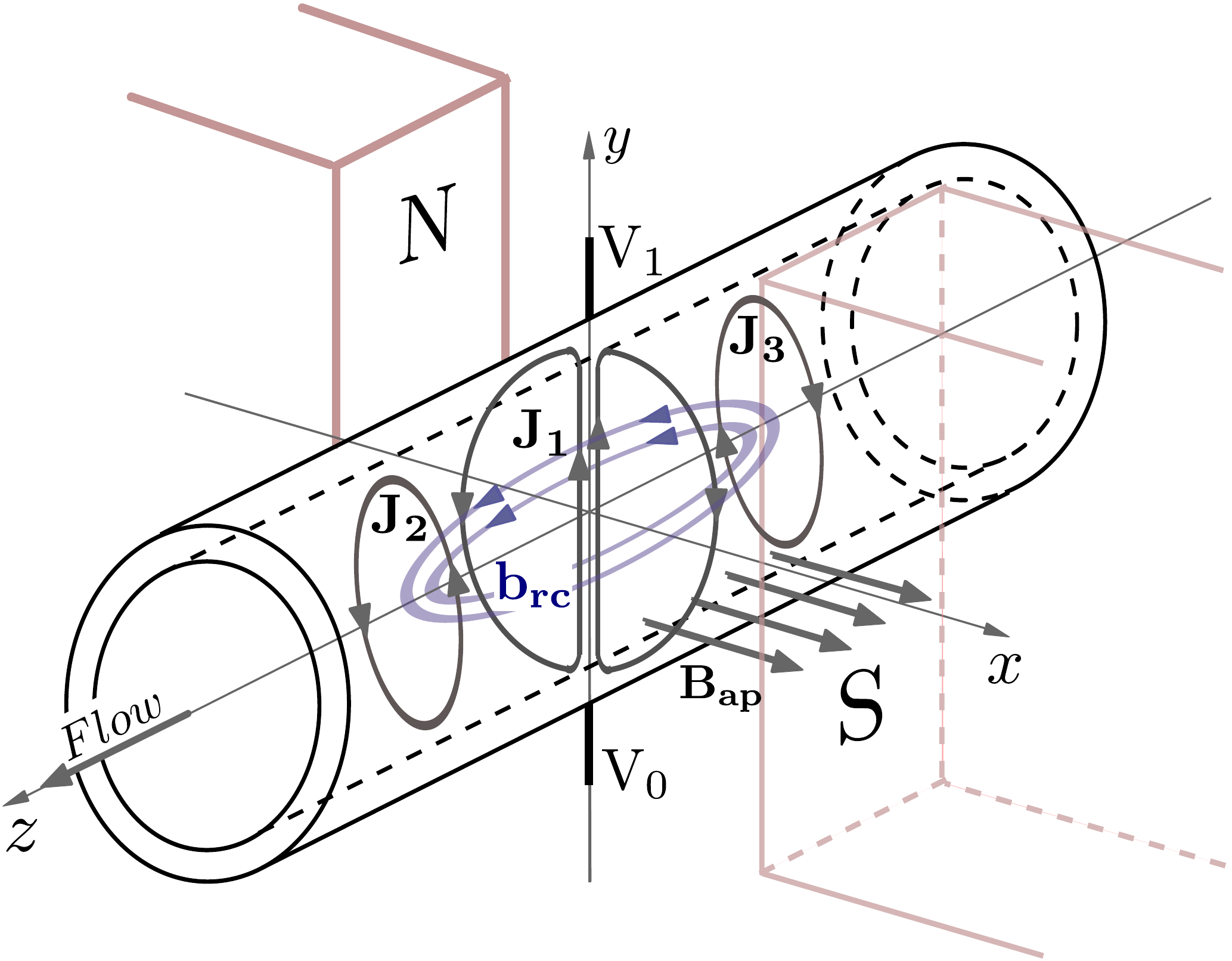}
		\caption{Schematic of Electromagnetic flowmeter. \cite{wikimedia1}}
		\label{em01}
	\end{figure}

	Fig. \ref{em01}. shows the schematic of  electromagnetic flowmeter. The measurement probes ($V_1 ~\& ~V_0$) are placed perpendicular to both magnetic field and flow direction. Circulating currents ($\bf{J_1, J_2, J_3}$) are due to the spatial variation in the induced electric field. The reaction magnetic field ($\bf{b_{rc}}$) produced by these currents cancels the applied magnetic field ($\bf{B_{ap}}$) in the upstream region and aids it in the downstream side. This cross-magnetizing effect apparently shifts the effective magnetic field along the flow direction. Generally in the analysis, the applied or the ambient and the reaction magnetic fields are separated.

The governing equations in terms of the magnetic vector potential $\bf{A}$ of the reaction field and the electric scalar potential $\phi$ arising out of current flow are given by \cite{emfmbook}, \cite{fmbase}:
\begin{equation} \label{eqge1}
\begin{split}
\nabla \cdot (\sigma \nabla \phi) - \nabla \cdot (\sigma ~ {\bf{u}} \times \nabla \times {\bf{A}}) =  \nabla \cdot (\sigma ~ \bf{u} \times \bf{B_{ap}})
\end{split}
\end{equation}	
\begin{equation} \label{eqge2}
\sigma \nabla \phi ~-~  \dfrac{1}{\mu} \nabla^2 {\bf{A}} - \sigma~ {\bf{u}} \times \nabla \times {\bf{A}} = \sigma~ {\bf{u}} \times {\bf{B_{ap}}}
\end{equation}	
where, $\mu$ is the magnetic permeability, $\sigma$ is the electrical conductivity and $\bf{u}$ is the velocity function of the fluid flow. {The  relative strength of the reaction magnetic field ($\bf{b_{rc}}$) to the applied magnetic field ($\bf{B_{ap}}$)  is indicated by the magnetic Reynolds number  $ R_m = \mu \sigma u_z D_h $ where, $D_h$ is the hydraulic diameter of the pipe }\cite{emfmbook}, \cite{fmbase}.

{For problems with $R_m < 1$, the reaction magnetic field ($\bf{b_{rc}}$) and hence the cross magnetizing effect becomes negligible. In such a situation, the last term on the Left Hand Side (LHS) of }(\ref{eqge1}) { becomes negligible and the resulting equation can be solved independently. It has been reported in [1] that whenever the length of the magnet is $>1.5$ times the pipe diameter, a two dimensional (2D) analysis  across the cross section perpendicular to the flow is permissible. For such cases the analytical solution can be found} \cite{fm2d2}, \cite{fm2d1}. In liquid metals, however, the conductivity is very high and hence the induced currents are large, which leads to strong cross-magnetizing effects. As a result, a full three dimensional (3D) analysis will be required. 
Due to the complexity in handling reaction field and the flow-geometry, numerical techniques like Galerkin Finite Element Method (GFEM) is generally employed. 

It is well known that GFEM, similar to the central-difference scheme, is averaging in nature and hence it is diffusive \cite{cdbook}.It is shown to become numerically unstable whenever the convection starts dominating over the diffusion. This numerical instability is widely addressed in the fluid dynamics literature dealing with transport equation \cite{cdbook}. Streamline Upwind Petrov Galerkin (SUPG) scheme \cite{supg1},  Galerkin Least Squares (GLS)  \cite{gls1}, Finite Increment Calculus (FIC) \cite{fic1} and Multiscale scheme \cite{multiscale1} are suggested for stabilizing the solution. 
{The same upwinding schemes have also been adopted for electromagnetic problems. For example the moving conductor problem has been analyzed in }\cite{mc1,mc2,mc3,mc4,mc5} {using the upwinding Petrov-Galerkin scheme.}

A similar numerical instability is also encountered with the GFEM simulation of the magnetic flowmeter. The SUPG scheme has been successfully employed in situations involving magnetic Reynolds number in the range 4 to 28 \cite{fmbase}. However, its performance for higher magnetic Reynolds number is yet to be quantified. Apart from this, SUPG involves more computation for higher order elements and also it is difficult to find the stabilization parameters for elements with order higher than quadratic \cite{quadp1}, \cite{quada1}. The present work basically aims to overcome these difficulties in the finite element simulation of electromagnetic flowmeters. 

{The investigations on the numerical instability, as well as, that for the remedial measures, were always carried out with the one dimensional version of the original problem} \cite{cdbook},  \cite{quadp1}, \cite{up1}. {  This is because, the required analysis in 2D or 3D is almost impractical and further, the one dimensional (1D) version of the fluid dynamics/thermal problems contained all the required features of the original problem.

 Following the same philosophy, the required analysis for the present work will be carried out using 1D version of the problem.} Both finite difference and Z-transform approaches will be employed for the analysis of numerical instability. From an insight thus obtained, a novel and stable scheme will be proposed for the classical GFEM. Subsequently the proposed method will be applied to the original flowmeter problem and the stability of the scheme will be demonstrated.


\section{\bf{Present work}}
For the theoretical investigation on the source of numerical instability, the analytical solution of the global set of equations given by GFEM is required. However, such an exercise is nearly impractical to be carried out in 2D or 3D and hence it is customary to resort to a 1D version of the problem \cite{cdbook},  \cite{quadp1}, \cite{up1}.

Following the same, a reduced one dimensional problem is considered. For this, conducting fluid is assumed to occupy the whole space with a spatially uniform velocity $u_z$ in the $z$-direction. The magnetic field is applied only in $x$-direction and is defined by,
\begin{equation}\label{bap1d}
B_x(z) = \left\{
        \begin{array}{ll}
            0 & \quad 0 \leq z < a \\
            B & \quad a \leq z \leq b \\
            0 & \quad b < z \leq L \\
        \end{array}
    \right.
\end{equation}
{where, $L$ is the length of the analysis domain and $B_x$ exists between $a $ and $ b $.}

 With these imposed conditions, the field variables cease to have any variation along $x$ and $y$ directions. The governing equations (\ref{eqge1}) \& (\ref{eqge2}) therefore reduces to, a single ordinary differential equation (ODE) in terms of $A_y$, the only non-vanishing component of $\bf{A}$:
\begin{equation}\label{eq1d1x}
~~-\dfrac{d^{2}A_{y}}{dz^2} + \mu \sigma u_z \dfrac{dA_y}{dz} = \mu\sigma u_z B_x ~
\end{equation}

Left hand side of the equation (\ref{eq1d1x}) has the same structure as the one used in fluid dynamics literature for investigating the numerical instability \cite{cdbook}, \cite{supg1}, \cite{fic1}.


\subsection{Analysis on instability}
	
	Application of the Galerkin finite element method (GFEM/FEM) and the finite  difference method (FDM) to (\ref{eq1d1x}), both results in same set of difference equation \cite{cdbook}, \cite{eddy1d1}.
	
	The difference equation for $n^{th}$ node,
	\begin{equation}\label{eq1dd2}
	 (-1-Pe)A_{y(n-1)} ~+~ 2 A_{y(n)} ~+~ (-1+Pe)A_{y(n+1)} \\~=~ 2Pe\Delta z  ~B_{x(n)}
	\end{equation}
{	where the Peclet number, $Pe = (\mu \sigma u_z \Delta z)/2$ and $\Delta z$ represents the element length. It may be worth noting here that, the oscillatory property of the numerical solution is basically a function of $Pe$ rather than its constituent parameters taken in isolation. Hence the allied literature considers $Pe$ as the independent variable in the required study }\cite{cdbook}, \cite{quadp1}, \cite{quada1}. When $Pe>1$, a root of the above difference equation becomes negative and it has been identified to be the source of numerical oscillation \cite{debook1}, \cite{dcbook2}.
	
	 The instability problem can also be analyzed by bringing the tools from the control system theory. The difference equation is transformed to the frequency domain for an easier analysis.
	
	The z-transform of (\ref{eq1dd2}),
	\begin{equation}\label{eq1db1z}
	 \Big((-1-Pe)Z^{-1} + 2 + (-1+Pe)Z~\Big) A_y = 2Pe\Delta zB_x
	\end{equation}
	(\ref{eq1db1z}) can be written in transfer function form,
	\[ \dfrac{A_y}{B_x} ~=~ \dfrac{2Pe\Delta z~}{-1+Pe}~~ \dfrac{Z}{Z^2~+~\dfrac{2}{-1+Pe}Z~+~\dfrac{-1-Pe}{-1+Pe}} \]

	\begin{equation}\label{bxzfinal}
	\dfrac{A_y}{B_x} ~=~ \dfrac{2Pe\Delta z}{-1+Pe}~~ \dfrac{Z}{(Z~-~1)~(Z~-~\dfrac{-1-Pe}{-1+Pe})} 
	\end{equation}
	when $Pe >> 1$
	\[ \dfrac{A_y}{B_x} ~\approx~ \dfrac{2 Pe\Delta z}{-1+Pe}~~ \dfrac{Z}{(Z~-~1)~(Z~+~1)} \]
	
	The above transfer function has poles at -1 \& +1. The pole located at '-1' is responsible for numerical oscillations \cite{dcbook2}. This observation is not specific to any particular excitation.
	
	In control systems, the controller design is always coupled with the pole-zero cancellation. Zeros of the transfer function  (\ref{bxzfinal}) can be modified by changing the input to a suitable form. A novel way to bring the desired zeros of the transfer function is to express the input magnetic field $B_x$ in terms of the associated magnetic vector potential $A_{sy}$.
	
	Substituting $B_x = -dA_{sy}/dz$ in the equation (\ref{eq1d1x}), in the difference form:
	\begin{equation}\label{eq:ay_diff}
	 (-1-Pe)A_{y(n-1)} ~+~ 2 A_{y(n)} ~+~ (-1+Pe)A_{y(n+1)} \\ ~=~  Pe(A_{sy(n-1)} ~-~ A_{sy(n+1)}) 
	\end{equation}

	The corresponding transfer function is:
	\begin{equation}\label{eq1da1}
	\dfrac{A_y}{A_{sy}} ~=~ \dfrac{-Pe}{-1+Pe} \dfrac{(Z~-~1 )~(Z~+~1 )}{(Z~-~1)~(Z~-~\dfrac{-1-Pe}{-1+Pe})} 
	\end{equation}
	
	when $Pe >> 1$
	\begin{equation}\label{eq1da2}
	 \dfrac{A_y}{A_{sy}} ~\approx~ \dfrac{-Pe}{~Pe} ~\dfrac{(Z~-~1 )~(Z~+~1 )}{(Z~-~1)~(Z~+~1)}
	\end{equation}
	Or
	\begin{equation}\label{eq1da3}
	\dfrac{A_y}{A_{sy}}~ \approx~ -1 
	\end{equation}
	
	A perfect pole-zero cancellation occurs for $Pe >> 1$, which ensures absolute stability. For very low Peclet numbers ($<1$), the intrinsic accuracy of GFEM is left unaltered. 
However, for Peclet numbers in the range $1 - 30$, the pole-zero cancellation is not perfect and hence some oscillation can arise. In order to find the maximum amplitude of this oscillation and the corresponding value of $Pe$, further work is carried out.
\subsection{Boundary conditions for the one dimensional version of the problem}
\begin{figure}
		\centering
    	  	  \mbox{\subfloat[]{\label{ip1d1} \includegraphics[scale=0.7]{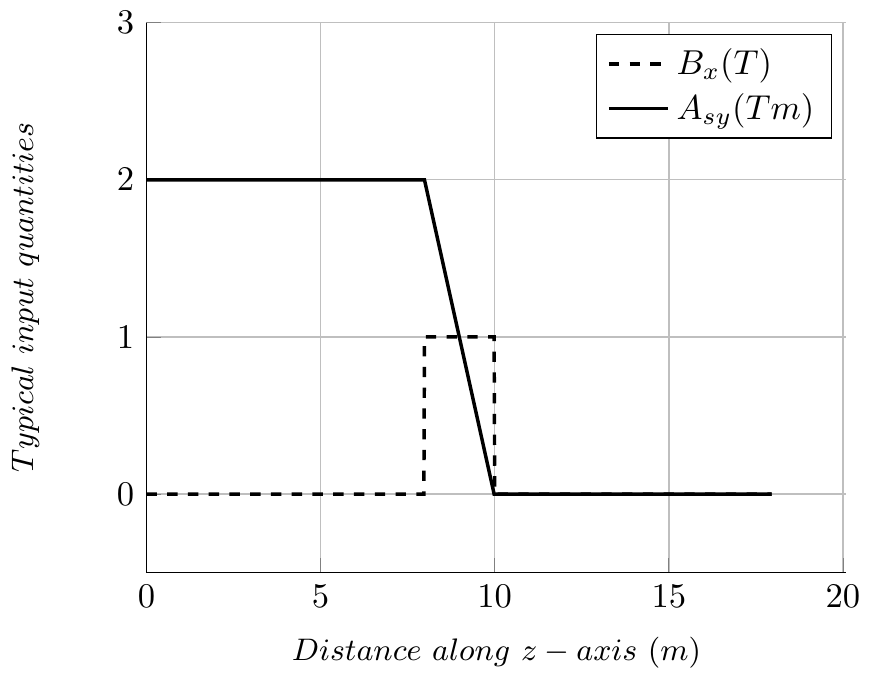}}}
    	  	  \mbox{\subfloat[]{\label{tdem1} \includegraphics[scale=0.3]{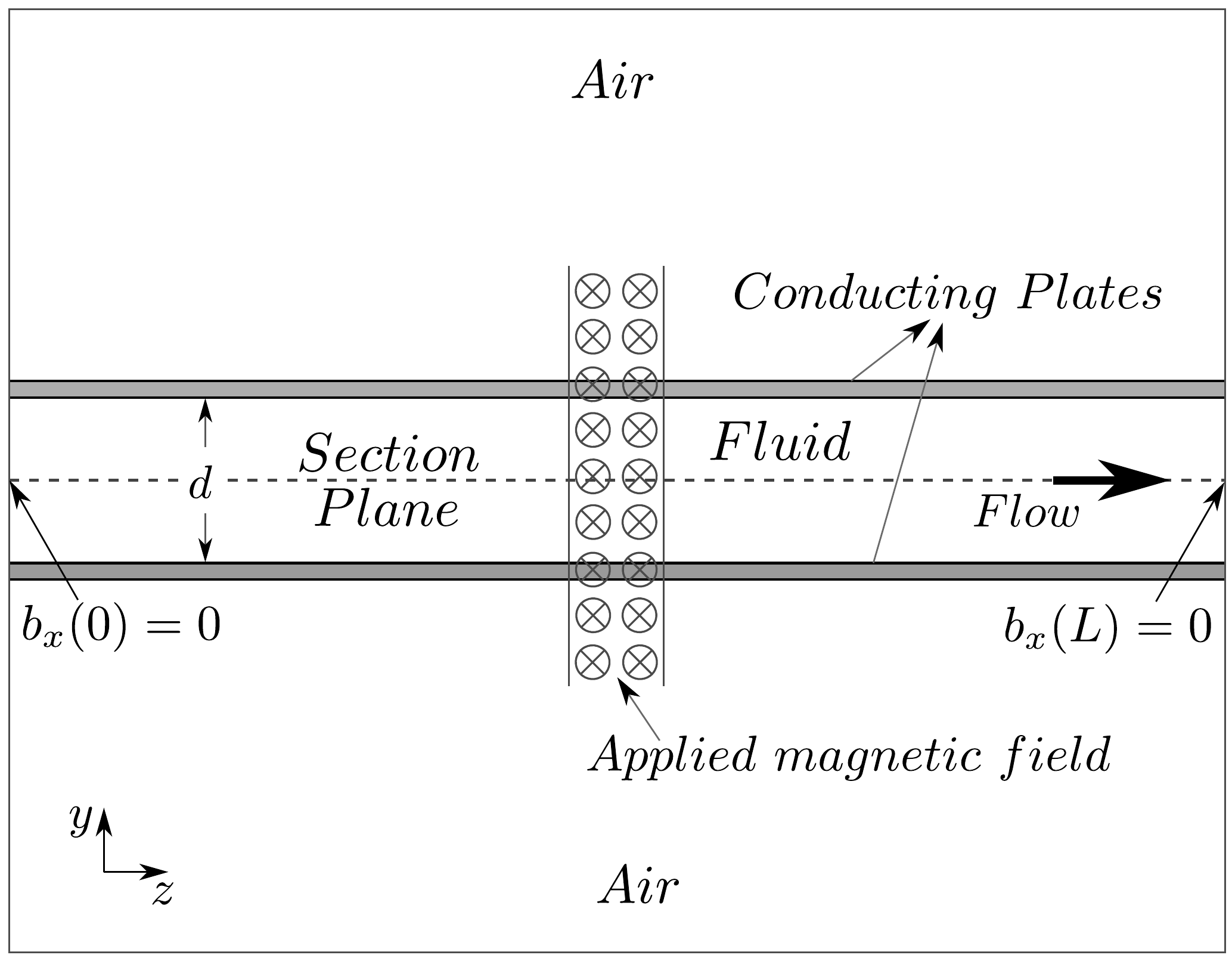}}}
		  \mbox{\subfloat[]{\label{tdbx1} \includegraphics[scale=0.8]{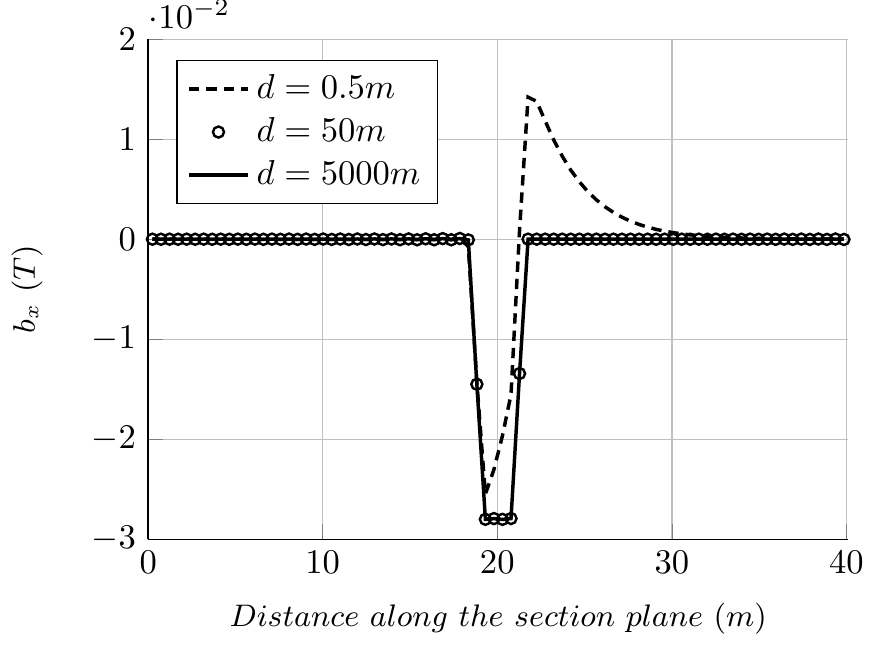}}}
		  \mbox{\subfloat[]{\label{tday1} \includegraphics[scale=0.8]{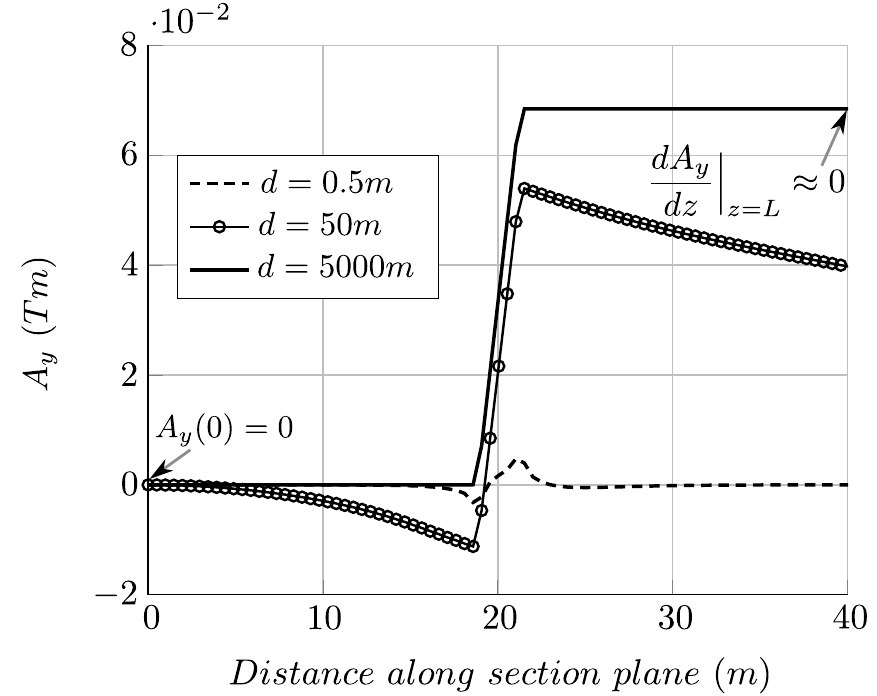}}}
		\caption{(a) Input quantity and its derivative form. (b) Schematic of 2D problem. (c) 2D: Computed reaction magnetic field for different plate separation distances $(d)$ with $u_z= 10 ~ms^{-1}$. (d) 2D: Computed magnetic vector potential for the same.}
		\label{tdrs}
	\end{figure}

The analytical and numerical solutions of difference equations (\ref{eq1dd2}) \& (\ref{eq:ay_diff}) will be considered. The input parameters employed are $\mu = 4 \pi \times 10^{-7} ~Hm^{-1}$, $\sigma = 7.21 \times 10^{6} ~Sm^{-1}$, $u_z$ and discretisation length ($\Delta z$) are kept as variables to make a  study on the oscillations. The applied magnetic field $B_x$ is presented in fig. \ref{ip1d1} along with the corresponding vector potential $A_{sy}$.
%

Because of the infinite spatial extension of the current in the one dimensional version of the problem, it cannot be directly mapped on to the original problem. Therefore, it becomes necessary to arrive at the appropriate boundary conditions for the same.


For this purpose, a two dimensional version of the problem is considered, which deemed to have adequate representation of the original problem. As shown in the fig. \ref{tdem1}, except for the fact that the fluid is now confined to the region between two finitely conducting plates kept parallel to the $xz$ plane, all other aspects are similar to 1D problem. The velocity of the fluid is assumed to be uniform in the gap. In order to apply the correct boundary conditions for the reaction magnetic field, the air region outside the plates up to a distance of $5d$ in $y$ direction  (where, $d$ is the separation distance between the plates) is included in the analysis. { The axial length of the computational domain is set to $10d$. For the FEM discretisation linear quadrilateral elements are employed. The element size along the  axial direction is kept constant while that along the $y$ direction is varied to accurately capture the current flow near the pipe wall. The element size selected always satisfied $Pe<1$.}

Even though the spatial extension of the current in the 2D problem is also infinite in the $x$ direction, the induced currents form closed loops in the $yz$ plane. Hence, the boundary condition ($\bf{A}$ = 0) which is relevant to the original problem is enforced at far distances of the upstream and downstream sides (fig. \ref{tdem1}). For a quantitative assessment of the reaction magnetic field and induced currents, numerical solution of the 2D problem is sought.

Selected results from the simulation are presented in {fig.} \ref{tdbx1} \& \ref{tday1} {for different gaps ($d$) between the plates.} In line with the expectation, it was found that for a wide range of velocities and the separation gap between the plates, the total current crossing the bisecting plane (shown in fig. (\ref{tdem1})) is zero. {This ensured that the reaction magnetic field $b_{x}$ vanishes at far distances along the downstream side. 
Further the increase in the gap ($d$) between the plates, $b_x$ is almost confined close to the region where the applied magnetic field exists.}


The solution for the 2D version of the problem has clearly shown that magnetic field is dragged only along the downstream side, while it is always confined at the upstream end (for the range of velocities of interest, which correspond to $\mu \sigma u_z > 10$). Based on this for the 1D problem, $A_y$ is specified to be zero at the upstream boundary and the vanishing $b_x$ leads to  $dA_y/dz = 0$ at the downstream boundary. 

\begin{figure}
		\centering
		  \mbox{\subfloat[]{\label{fdm1d1:3000} \includegraphics[scale=1]{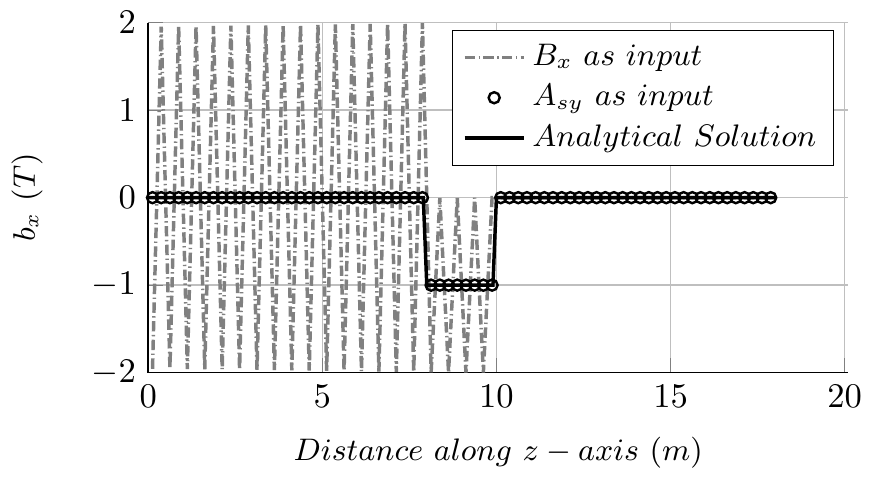}}}
		  \mbox{\subfloat[]{\label{fdm1d1:0003} \includegraphics[scale=1]{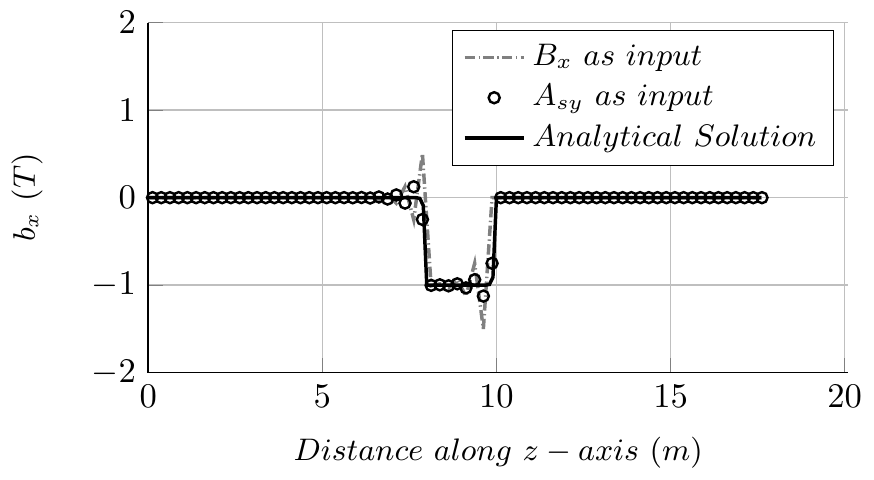}}}
		\caption{FDM/FEM solution.  $~$ (a) Pe = 3000, {$\Delta z = 0.25$.}$~~$ (b) Pe = 3, {$\Delta z = 0.25$.}}
		\label{fdm1d1}
	\end{figure}

Along with these boundary conditions and the input magnetic field defined in (\ref{bap1d}), the analytical solution of the ODE (\ref{eq1d1x}) is:
\begin{equation}\label{odeay}
A_y(z) = \left\{
        \begin{array}{ll}
            \dfrac{B}{k}(e^{-kb} - e^{-k(b-z)} \dots \\ ~~~~~~~~~- e^{-ka} + e^{-k(a-z)}) & \quad 0 \leq z < a \\
            \dfrac{B}{k}(1-e^{-ka} + e^{-kb} \dots \\ ~~~~~~~~~- e^{-k(b-z)}) + B(z-a) & \quad a \leq z \leq b \\
            \dfrac{B}{k}(e^{-kb} - e^{-ka}) + B(b-a) & \quad b < z \leq L \\
        \end{array}
    \right.
\end{equation}

where, $k = \mu \sigma u_z;$  $B =$ value of $B_x$ in $a \leq z \leq b$ as mentioned in (\ref{bap1d}). This can serve as the reference for the evaluation of the error in the numerical solution of the governing equation. Sample results obtained from the FEM are presented in  fig. \ref{fdm1d1} along with the analytical solution of the governing equation,  where the reaction magnetic field $b_x$ is calculated using the forward difference scheme:
\[ b_x(n) = -\dfrac{A_y(n+1)-A_y(n)}{\Delta z} \]
The solution is stable for $Pe>>1$ (practically verified for $Pe$ ranging from 30 to 30000) with $A_{sy}$ as input.  However, as mentioned earlier, in the mid-range $1<Pe<30$ oscillation exists. Sample results are presented in fig. \ref{fdm1d1}.

In order to identify the location (the value of $Pe$) and peak amplitude of the oscillation,  analytical solution of the FEM equation, which appears in the difference equation form, is performed in the next section.

\subsection{Location and value of the maximum error}
The governing difference equations can be rewritten as:

For magnetic field $B_x$ as the input:	
\begin{equation} \label{deqbxx}
rA_y(n-1)+(-1-r)A_y(n) +A_y(n+1)\\ = 2Pe\Delta z  ~B_{x}(n)
\end{equation}
Similarly for vector potential $A_{sy}$ as the input,	
\begin{equation} \label{deqasy}
rA_y(n-1)+(-1-r)A_y(n) +A_y(n+1)\\ = \frac{1-r}{2}(A_{sy}(n-1)-A_{sy}(n+1))
\end{equation}
where, $ r = (-1-Pe)/(-1+Pe)~$;
\medskip

Using the method described in  \cite{debook1}, the general solution of the equations (\ref{deqbxx}) and (\ref{deqasy}) can be found as,
\begin{equation} \label{gsbxx}
A_y(n)~=~k_1~+~k_2~r^{n}~+~y_p(n)
\end{equation}
where, $y_p(n)$ is the particular solution and 
$~~ k_1,~k_2$ are parameters of the complementary solution. 

{The direct solution of the difference equation with piecewise input is difficult and hence the problem  domain is divided into five sub-domains ($B,~C,~D,~F ~\&~G$) as shown in figs}. \ref{ip1d_bx} \& \ref{ip1d_asy} {of Appendix.} The solution steps are also detailed in the Appendix along with its validation.

The oscillation in the numerical result can be found both analytically and numerically to occur at the end of domain $C$. Based on this, for finding the peak amplitude of the oscillation and corresponding $Pe$, solution for domain $C$ i.e., equation (\ref{gsc2}) is considered. The reaction magnetic field at the end of domain $C$:
\begin{equation}\label{oscan1}
 b_{a} = \frac{A_y(m_c -1) - A_y(m_c)}{\Delta z}\\
  = \frac{c_2(r^{m_c-1}-r^{m_c})}{\Delta z}- \frac{\lambda}{\Delta z}\\
\end{equation}

Equation (\ref{oscan1}), has a constant part given by $\lambda\ \Delta z$ and an oscillatory part given by $ (c_2(r^{m_c-1}-r^{m_c}))/ \Delta z$.

Comparing with the analytical solution given in equation (\ref{odeay}), oscillatory part can be identified to be the major part of the error in the numerical solution. For the range of velocities of interest (i.e. $\mu \sigma u_z > 10$) this part amounts to almost the total error.

Let,
\begin{equation}\label{oscan2}
\widehat{b_{a}} =\frac{c_2(r^{m_c-1}-r^{m_c})}{\Delta z}\\
\end{equation}
	substituting (\ref{c2a}) in (\ref{oscan2}),
\begin{equation}\label{oscan3}
\widehat{b_{a}} ~=~ B \frac{1+r}{2r^2} ~=~  \frac{B(1-Pe)}{(1+Pe)^2} \\
\end{equation}	
Similarly, for $B_x$ input,
\begin{equation}\label{oscan3b}
\widehat{b_{b}} ~=~ \frac{B}{r}~=~\frac{B(1-Pe)}{1+Pe}\\
\end{equation}
	\begin{figure}[htp]
		\centering
\includegraphics[scale=1]{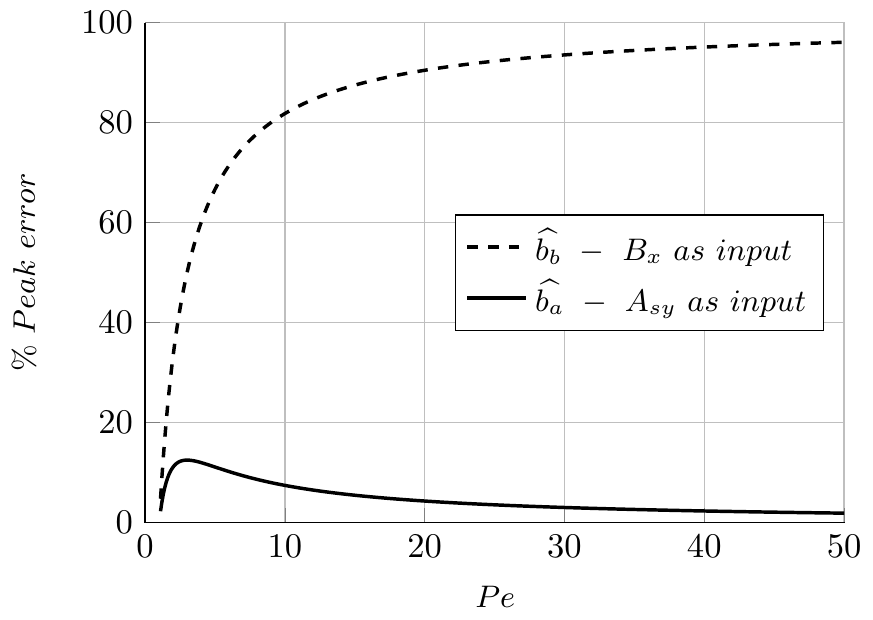} 
		\caption{\% peak error in the numerical solution for a wide range of  Peclet numbers.  }
		\label{babb1}
	\end{figure}	

The magnitude of the oscillation present in the analytical solution of the difference equation is given by (\ref{oscan3}) \& (\ref{oscan3b}) for the respective inputs (i.e., $A_{sy}$ \& $B_x$ respectively). The error obtained from the above is plotted in fig. \ref{babb1}. The peak error, with input field specified in terms of the vector potential occurs at $ Pe = 3$ and its magnitude in \% is $1/8$.

The above numerical exercise has once again confirmed that the proposed scheme is very stable for large flow rates. Also, even in the midrange of flowrates (and hence  Pe), the error in the numerical results is lower than that with GFEM and further it can be controlled by opting for different discretisation.


\subsection{Performance with quadratic elements}
\begin{figure}
		\centering
		  \mbox{\subfloat[]{\label{veraip1:a} \includegraphics[scale=1]{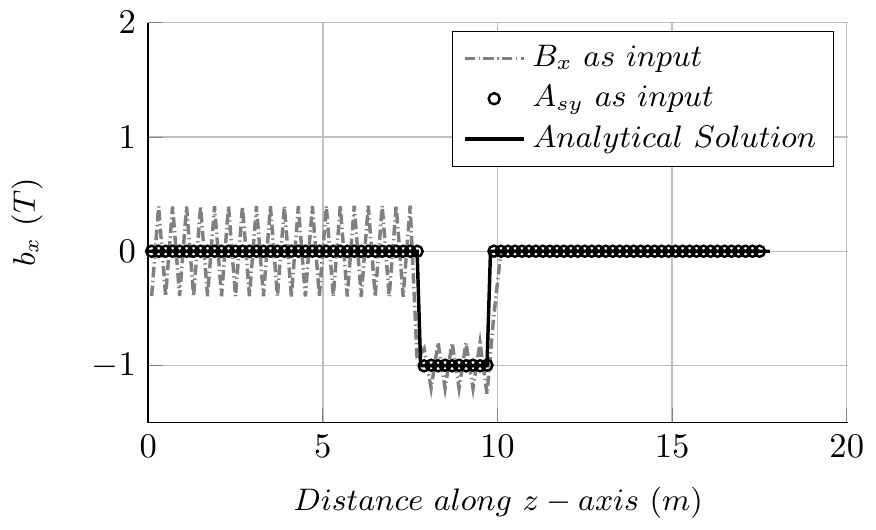}}}
		  \mbox{\subfloat[]{\label{veraip1:b} \includegraphics[scale=1]{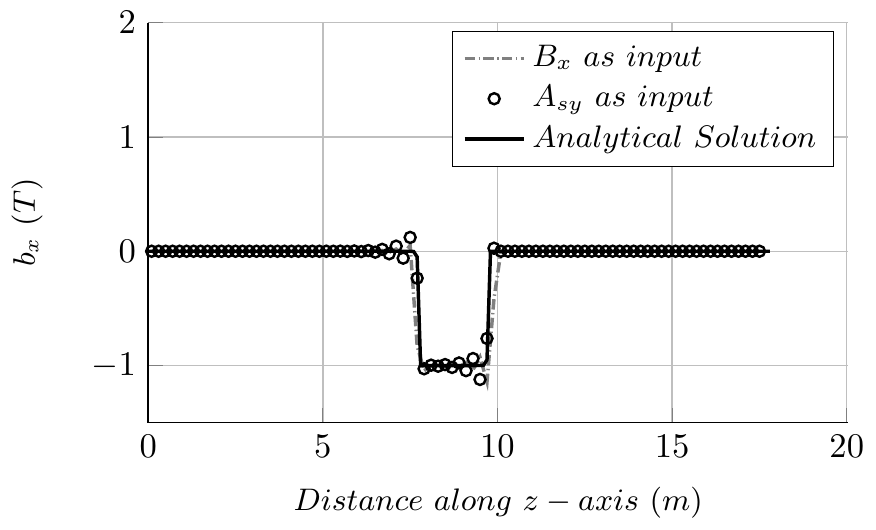}}}
		\caption{Verification for quadratic elements.  $~$ (a) Pe = 3000.$~~$ (b) Pe = 3.}
		\label{fig_qd}
	\end{figure}

Unlike that with the first order elements, SUPG scheme for higher order elements requires more computation. On the other hand, the proposed scheme is free of such issues and when implemented for second order elements, its performance is found to be equally good. Sample numerical results are presented in fig. \ref{fig_qd} as an evidence for the same. 
\subsection{Analysis for the flowmeter}

\begin{figure}
		\centering
		  \mbox{\subfloat[]{\label{ip3d_bx} \includegraphics[scale=0.8]{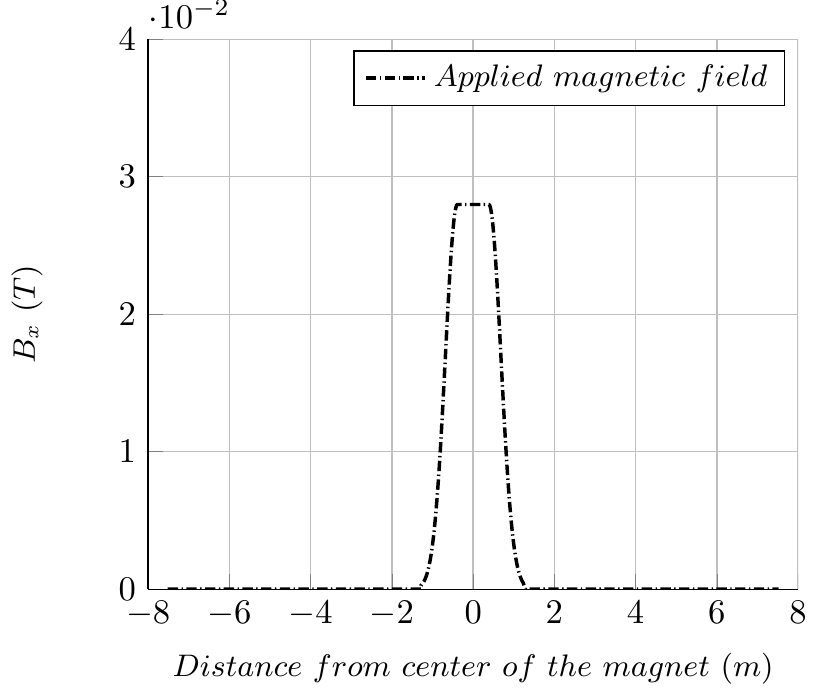} }}
		  \mbox{\subfloat[]{\label{pipeaxisplot} \includegraphics[scale=0.8]{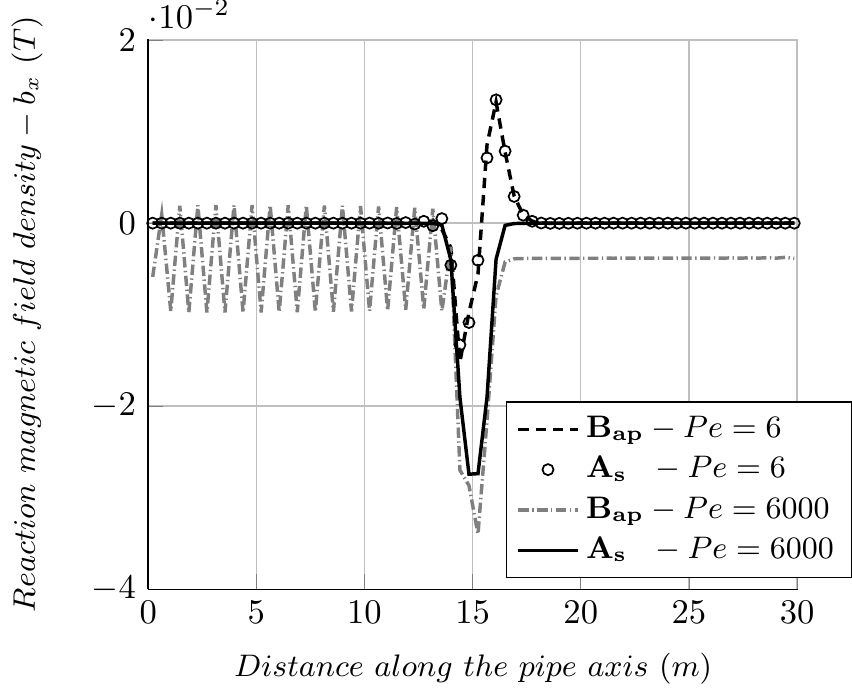}  }}
		  \mbox{\subfloat[]{\label{oscbx1} \includegraphics[scale=0.5]{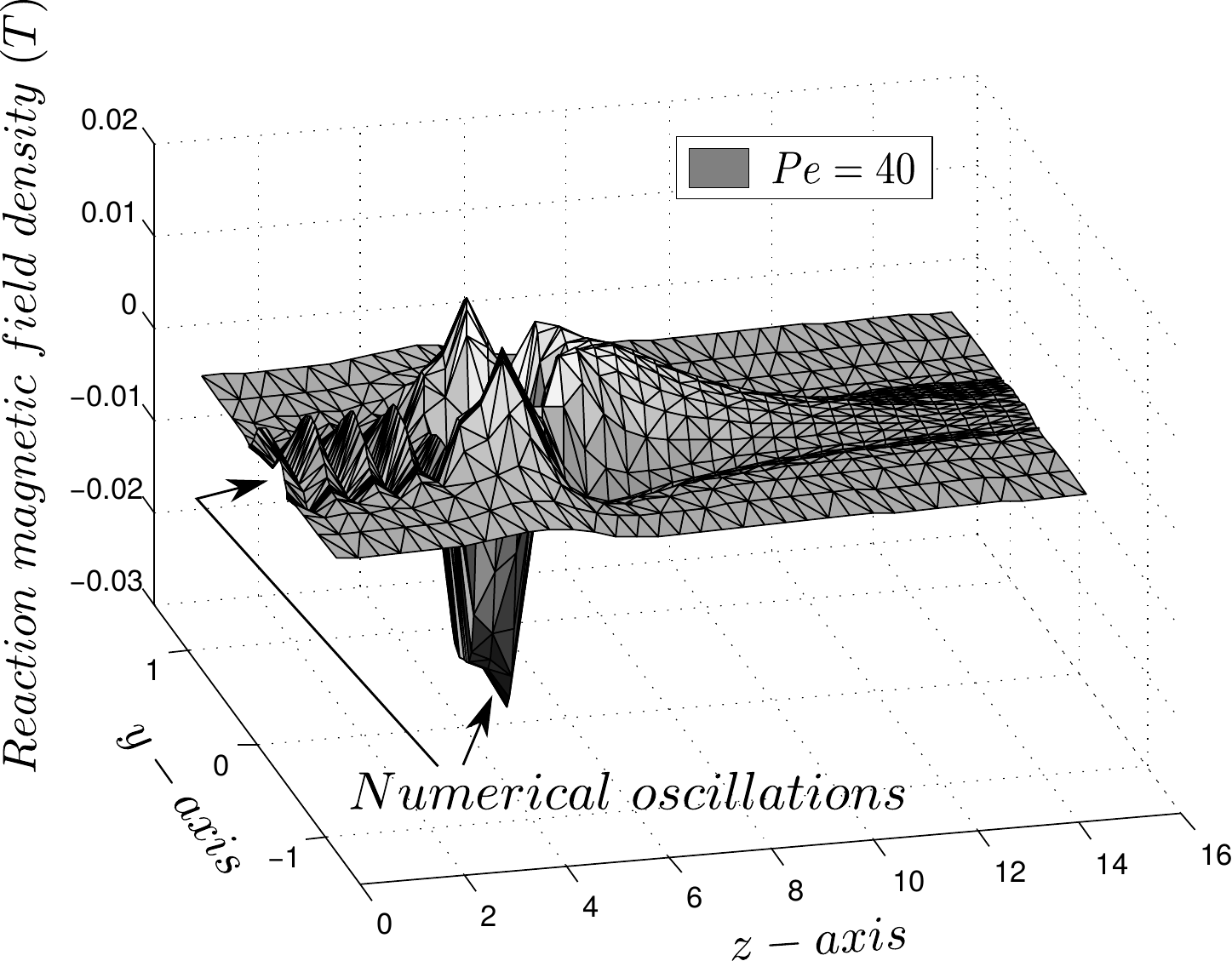}}}
		  \mbox{\subfloat[]{\label{stbbx1} \includegraphics[scale=0.5]{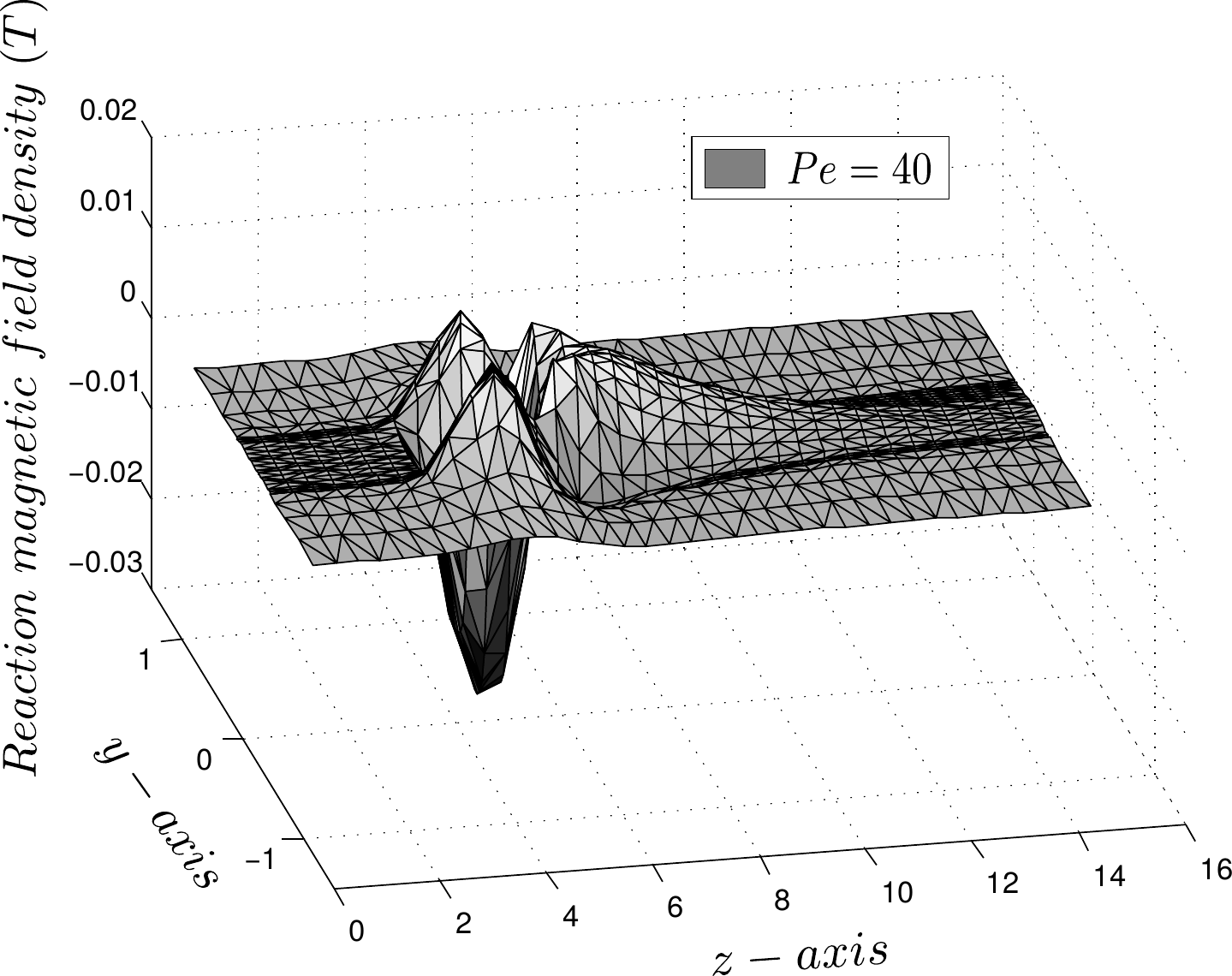}}}
		\caption{(a) Axial distribution of $B_x.~$ (b) $x$-component of ${\bf{b_{rc}}}$ along the pipe axis ($z$-axis) for different Peclet numbers with ${\bf{B_{ap}}}$ \&  $\nabla \times {\bf{A_s}}$  as inputs.$~$ (c) $x$-component of reaction magnetic field ${\bf{b_{rc}}}$ in  $~x$ = 0 plane with ${\bf{B_{ap}}}$ as input. $~$ (d) $x$-component of reaction magnetic field ${\bf{b_{rc}}}$ in  $~x$ = 0 plane with $ \nabla \times {\bf{A_s}}$  as input.}
		\label{stb3d}
	\end{figure}

Up till now all the analysis was limited to one dimensional version of the problem and therefore it was deemed necessary to scrutinize the proposed scheme with the original problem. For this, governing equations (\ref{eqge1}) \& (\ref{eqge2}) are solved using GFEM \cite{fmbase}. { The input magnetic field $\bf{B_{ap}}$ is set to have only the $x$-component $B_x$, whose magnitude varies only along the flow direction (as shown in fig} \ref{ip3d_bx}{). The corresponding vector potential $\bf{A_s}$ will have only $y$ component. It is readily available if one employs numerical simulation for the input magnetic field or it can be obtained from $-dA_{sy}/dz = B_x$ when the measured flux density is employed.}

 {For the numerical experiment, a non-magnetic stainless steel pipe with inner and outer diameters $ 517~mm ~\&~ 557 ~mm $ respectively is considered. The conductivity of the pipe material is $ 1.16 \times 10^{6} ~Sm^{-1}$ and that of the liquid sodium inside is $ 7.21 \times 10^{6} ~Sm^{-1}$. It is true that due to cavitation and other associated problems, velocities beyond few to few tens of meters  are impractical with liquid metals. Nevertheless,  simulations are carried out for velocities ranging up to $3000~ ms^{-1}$ (which corresponds to $R_m=14851$), solely to demonstrate the  robustness of the proposed scheme. The intention was to consider the possible application of the proposed scheme for general moving conductor problems. Incidentally, $R_m$ is related to the Peclet number through $ Pe = R_m \Delta z/2 D_h$.}

Sample simulation results are compared with that obtained from the original formulation in fig.\ref{stb3d}. The reaction magnetic field for $ Pe~=~40$ as given by two approaches is presented in fig. \ref{oscbx1} \& \ref{stbbx1},  while fig. \ref{pipeaxisplot} presents reaction magnetic field along the pipe axis for different values of $Pe$. {From the results obtained from extensive simulations, it is observed that for the  midrange of $Pe$ there exists a small amount of oscillation, which asymptotically vanishes with increase in $Pe$. These observations are in line with the inference drawn earlier from the 1D analysis. It must be noted here that all the analysis reported here considers only an axially varying applied magnetic field, which generally holds true for flowmeters.}

\section{\bf{Summary and Conclusion}}
Theoretical evaluation of the sensitivity of electromagnetic flowmeter seems to be the best possible choice especially when it is to be used for the  measurement of liquid metal flows.
The commonly employed Galerkin finite element formulation is known to become  unstable for larger flow rates.  SUPG scheme is generally suggested in the pertinent literature for overcoming this problem. However SUPG scheme requires computation of the stabilization parameter, which involves a lot more calculations with higher order elements. In addition, it is difficult to arrive at the stabilization parameters for elements with order beyond quadratic  \cite{quadp1}, \cite{quada1}.

 By analyzing the one dimensional version of the problem a novel scheme has been devised which is free of above mentioned difficulties. In this scheme, classical GFEM is retained intact and only the input magnetic field is restated in terms of the associated vector potential. {The analytical solution of the FEM, for the proposed scheme indicate that solution exhibits a small amount of oscillation in the mid-range of flow rates ($1<Pe<30$) and this oscillation asymptotically vanishes with the increase in $Pe$ (for $Pe > 3$).  The maximum error due to this numerical oscillation has been quantified analytically. Even the maximum error, which is shown to occur at $Pe~ =~ 3$, is lower than what is found with flux density as the input. Finally, the proposed scheme is applied to the original flowmeter problem and the simulation results indicate that inferences drawn on the 1D version of the problem remain valid even for the 3D case. In summary, a simple and robust scheme has been proposed for the FEM solution of the flowmeter problem involving only an axially varying applied magnetic field.}


\appendix
\renewcommand*{\thesubsection}{\arabic{subsection}}
\section{Analytical solution of the difference equation}
\subsection{Analytical solution - $B_x$ as Input}

\begin{figure}
		\centering
		\mbox{\subfloat[]{\label{ip1d_bx} \includegraphics[scale=0.53]{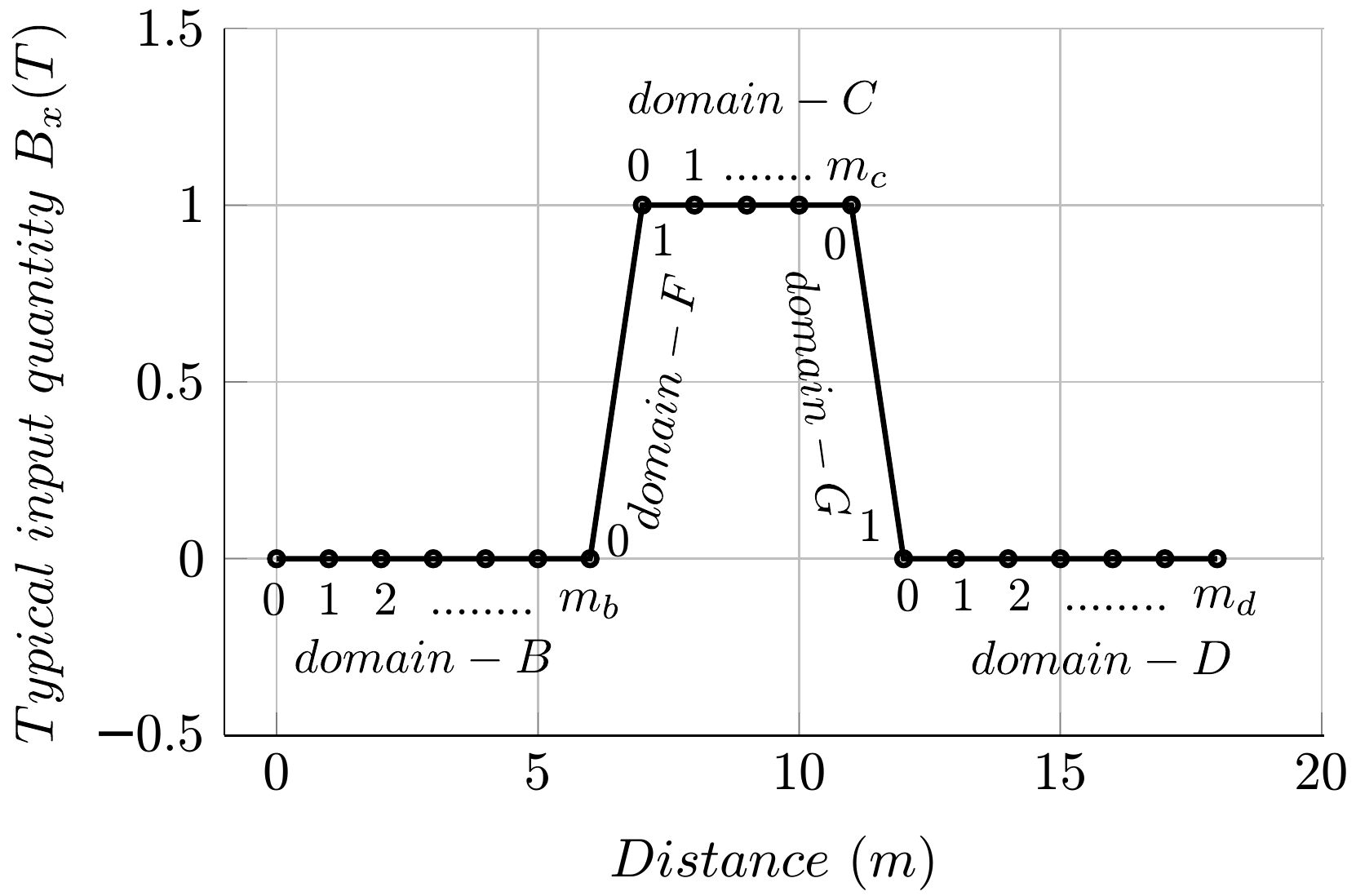}}}
		  \mbox{\subfloat[]{\label{verbip1:a} \includegraphics[scale=0.9]{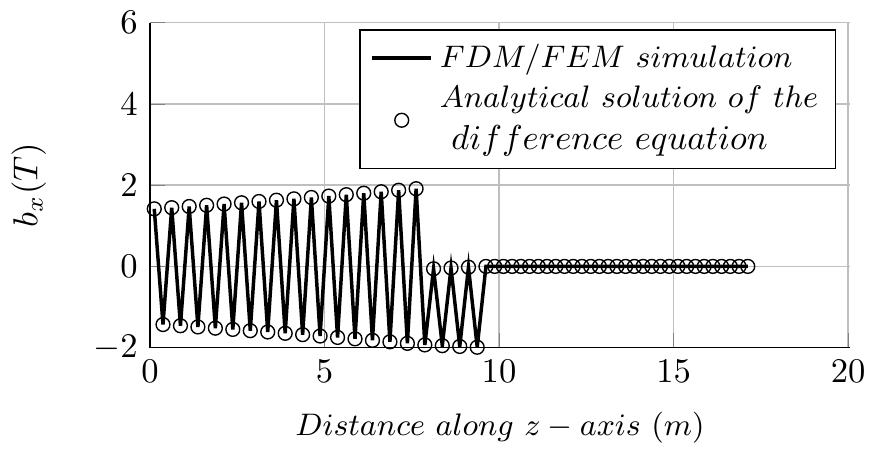}}}
		  \mbox{\subfloat[]{\label{verbip1:b} \includegraphics[scale=0.9]{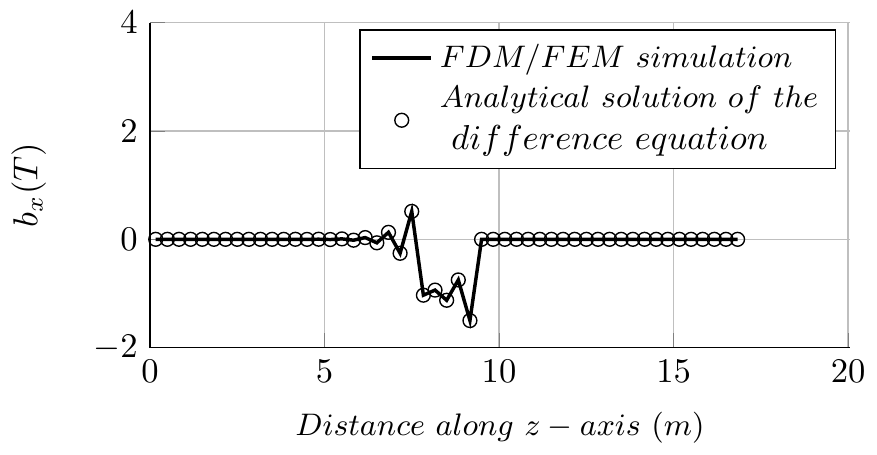}}}
		\caption{(a) Five sub-domains and their ranges $-$ $B_x$ as Input. (b) Validation of the analytical solution with  $B_{x}$ as input $~$ $Pe=200,~\Delta z = 0.25,~ m_c=6, ~m_b = 31, ~m_d = 31.~~$ (c) Validation of the analytical solution with  $B_{x}$ as input $~$ $Pe=3,~\Delta z = 0.33,~ m_c=4, ~m_b = 23, ~m_d = 23 $.}
		\label{verbip1}
	\end{figure}

The direct solution of a difference equation with piecewise-defined input is difficult. Hence the domain is divided into five sub-domains as shown in fig. {\ref{ip1d_bx}}

 Following notations $\lambda =B ~\Delta z,~r = (-1-Pe)/(-1+Pe),$ and $y = A_y$ are employed in the subsequent steps. In order to distinguish the solutions of each domain, the parameters of complementary solution are designated with their domain names ($b_1$, $b_2$,  $c_1$, $c_2$ $~$ etc.,). And the other variables such as $y$, particular solution $y_p$ and node number $n$ are suffixed with their corresponding domain names ($y_b$, $n_b$, $y_c$, $n_c$, $y_{pc}$ $~$ etc., ). The governing difference equations and their general solutions for different domains are listed below.
\medskip

For domain $B$ ($0 \leq n_b \leq m_b $):
\begin{eqnarray} \label{deqb1}
ry_b(n_b-1)+(-1-r)y_b(n_b) +y_b(n_b+1)=0 
\end{eqnarray}
\begin{equation} \label{gsb1}
\Rightarrow y_b(n_b)~=~b_1~+~b_2~r^{n_b}  
\end{equation}

For domain $F$ ($0 \leq n_f \leq 1 $):
\begin{equation} \label{deqf1}
ry_f(n_f-1)+(-1-r)y_f(n_f)+y_f(n_f+1)=~~~~\\ \lambda (1-r)n_f
\end{equation}
\begin{equation} \label{gsf1}
\Rightarrow y_f(n_f)~=~f_1~+~f_2~r^{n_f}~+~y_{pf}(n_f)
\end{equation}

For domain $C$  ($0 \leq n_c \leq m_c $):
\begin{equation} \label{deqc1}
ry_c(n_c-1)+(-1-r)y_c(n_c) +y_c(n+1)=\lambda (1-r) 
\end{equation}
\begin{equation} \label{gsc1}
\Rightarrow y_c(n_c)~=~c_1~+~c_2~r^{n_c}~+~y_{pc}(n_c)
\end{equation}

For domain $G$  ($0 \leq n_g \leq 1 $):
\begin{equation} \label{deqg1}
ry_g(n_g-1)+(-1-r)y_g(n_g)+y_g(n_g+1)=~~~~\\ \lambda (1-r)(1-n_g)
\end{equation}
\begin{equation} \label{gsg1}
\Rightarrow y_g(n_g)~=~g_1~+~g_2~r^{n_g}~+~y_{pg}(n_g)
\end{equation}

For domain $D$  ($0 \leq n_d \leq m_d $):
\begin{equation} \label{deqd1}
ry_d(n_d-1)+(-1-r)y_d(n_d) +y_d(n_d+1)=0
\end{equation}
\begin{equation} \label{gsd1}
\Rightarrow y_d(n_d)~=~d_1~+~d_2~r^{n_d}
\end{equation}
The difference equations in domains $F,C~ \&~ G$ are inhomogeneous and their particular solutions are computed \cite{debook1}. The complete general solutions of domains $F,C~ \&~ G$ are,

\begin{equation} \label{gsf2}
y_f(n_f)=f_1+f_2r^{n_f}+\frac{\lambda}{2} \frac{(r+1)}{(r-1)} n_f+
\frac{\lambda}{2} n_f^2
\end{equation}
\begin{equation} \label{gsc2}
y_c(n_c)=c_1+c_2r^{n_c}+\lambda n_c
\end{equation}
\begin{equation} \label{gsg2}
y_g(n_g)=g_1+g_2r^{n_g}+ \Big( 1 - \frac{r+1}{2(r-1)} \Big) \lambda n_g-
\frac{\lambda}{2} n_g^2
\end{equation}

	The global boundary conditions imposed on the extreme boundaries. It helps in eliminating two variables.
\begin{equation}\label{b1}
 y_b(0) = 0~~~\Rightarrow~~ b_1 = -b_2\\
\end{equation}
\begin{equation}\label{d2}
 y_d(m_d) = y_d(m_d+1)~~~\Rightarrow~~ d_2 = 0\\
\end{equation}

The intermittent boundaries should satisfy the equality condition and the governing equation itself at the joining nodes. At junction $G-D$,
	\begin{itemize}[label={}]
	\item Equality condition: 
		\begin{equation}\label{gdeq1}
		y_g(1) = y_d(0) ~~\Rightarrow 	~	g_1 + g_2 r+ \frac{\lambda}{1-r}  = d_1
		\end{equation}
	\item Satisfying, governing equation
		\begin{equation}\label{gdeq2}
		ry_g(0) - (1+r)y_d(0) + y_d(1) = 0~~ \Rightarrow
		g_1 + g_2   = d_1
		\end{equation}
	\end{itemize}

The above mentioned conditions are employed for other three junctions and six more equations are obtained. These six equations are solved along with equations (\ref{b1}), (\ref{d2}), (\ref{gdeq1}) \& (\ref{gdeq2}). The parameters of complementary solutions are obtained as:
\begin{flushleft}
\begin{eqnarray}
g_2 &=& \frac{\lambda}{(1-r)^2} \label{g2}\\
c_2 &=& \frac{\lambda}{r^{m_c}(1-r)}\label{c2}\\
f_2 &=& -\lambda \frac{1}{(1-r)^2} + \frac{c_2}{r}\label{f2}\\
b_2 &=& \frac{\lambda}{r^{m_b-1}(1-r)^2} + \frac{f_2}{r^{m_b}}\label{b2}\\
f_1 &=& b_1 + b_2 r^{m_b} - f_2\label{f1}\\
c_1 &=& f_1 + f_2 r + \lambda \frac{r}{r-1} -c_2\label{c1}\\
g_1 &=& c_1 + c_2 r^{m_c} + \lambda m_c - g_2\label{g1}\\
d_1 &=& g_1 + g_2 r+ \frac{\lambda}{1-r} \label{d1}
\end{eqnarray}
\end{flushleft}
	
The analytical solution of the difference equation is compared with the numerical solution obtained from FDM and FEM in fig. \ref{verbip1:a} \& \ref{verbip1:b}.
\subsection{Analytical solution with derivative of $A_{sy}$ as Input}
\begin{figure}
		\centering
		  \mbox{\subfloat[]{\label{ip1d_asy} \includegraphics[scale=0.95]{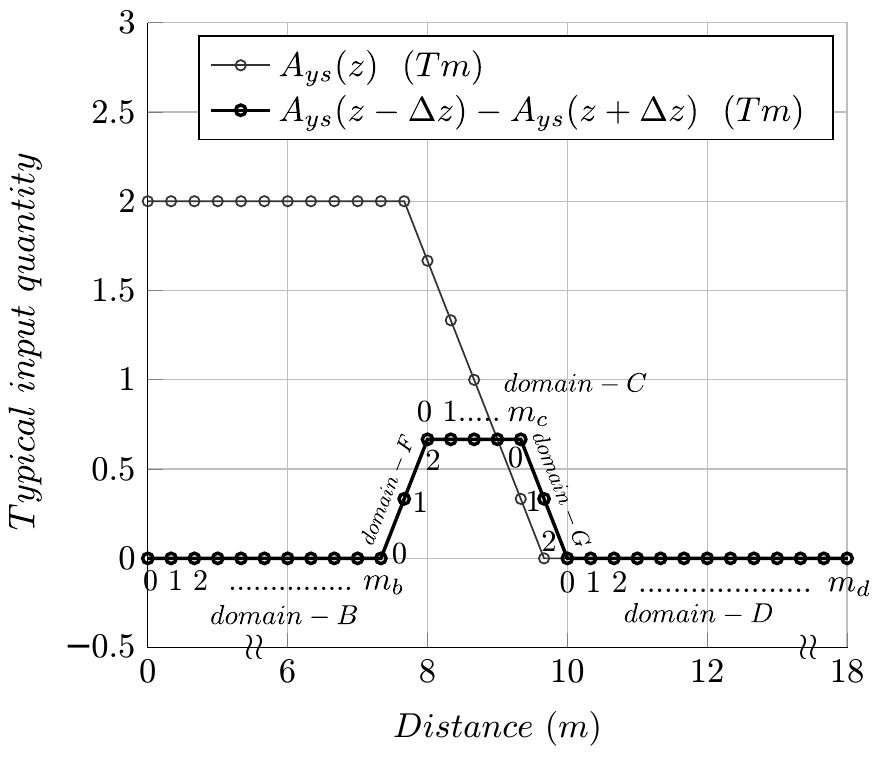}}}
		  \mbox{\subfloat[]{\label{veraip1:a} \includegraphics[scale=0.9]{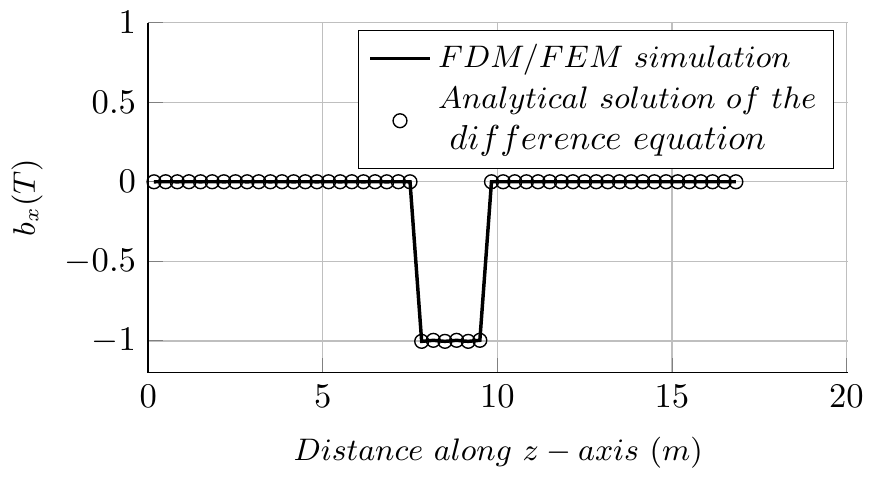}}}
		  \mbox{\subfloat[]{\label{veraip1:b} \includegraphics[scale=0.9]{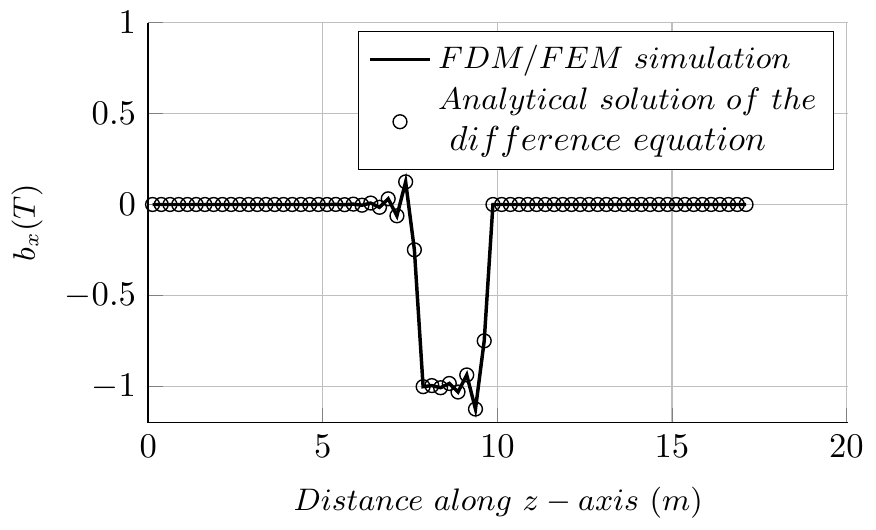}}}
		\caption{(a) Five sub-domains and their ranges $A_{sy}$ as input. $~$(b) Validation of the analytical solution with $A_{sy}$ as input $Pe=300,~\Delta z = 0.33,~ m_c=4, ~m_b = 22, ~m_d = 22$. $~$(c) Validation of the analytical solution with $A_{sy}$ as input $Pe=3,~\Delta z = 0.25,~ m_c=6, ~m_b = 30, ~m_d = 30$.}
		\label{veraip1}
	\end{figure}

The input in terms of $A_{sy}$ is plotted in fig. \ref{ip1d_asy}. Similar to $B_x$ input, the problem is divided into five sub-domains. 

The governing difference equations and their general solutions for domains $B,~C,~\&~D$ are same as that of the previous case, with $\lambda$ redefined as $ \lambda ~=~ (A_{sy}(n_c-1) - A_{sy}(n_c+1))/2$ and it is equivalent to $B~\Delta z $. $F~ \&~ G$ domain equations are listed as follows,
\medskip

For domain $F$  ($0 \leq n_f \leq 2 $):
\begin{equation} \label{deqf1a}
ry_f(n_f-1)+(-1-r)y_f(n_f)+y_f(n_f+1)=~~~~\\ \lambda  \frac{1-r}{2} n_f\\
\end{equation}
\begin{equation} \label{gsf2a}
\Rightarrow y_f(n_f)=f_1+f_2r^{n_f}+\frac{\lambda}{4} \frac{(r+1)}{(r-1)} n_f+
\frac{\lambda}{4} n_f^2
\end{equation}

For domain $G$  ($0 \leq n_g \leq 2 $):
\begin{equation} \label{deqg1a}
ry_g(n_g-1)+(-1-r)y_g(n_g)+y_g(n_g+1)=~~~~\\ \lambda (2-n_g) \frac{1-r}{2}
\end{equation}
\begin{equation} \label{gsg2a}
\Rightarrow y_g(n_g)=g_1+g_2r^{n_g}+ \Big( 1 - \frac{r+1}{4(r-1)} \Big) \lambda n_g-
\frac{\lambda}{4} n_g^2
\end{equation}

The computed values of parameters in the complementary solution are:
\begin{flushleft}
\begin{eqnarray}
g_2 &=& \frac{\lambda}{2r(1-r)^2} \label{g2a}\\
c_2 &=& \frac{\lambda(1+r)}{2r^{m_c+1}(1-r)}\label{c2a}\\
f_2 &=& \frac{-\lambda}{2r(1-r)^2} + \frac{c_2}{r^2}\label{f2a}\\
b_2 &=& \frac{\lambda r}{2r^{m_b}(1-r)^2} + \frac{f_2}{r^{m_b}}\label{b2a}\\
f_1 &=& b_1 + b_2 r^{m_b} - f_2\label{f1a}\\
c_1 &=& f_1 + f_2 r + \frac{\lambda(3r-2)}{2(r-1)} - \frac{c_2}{r}\label{c1a}\\
g_1 &=& c_1 + c_2 r^{m_c} + \lambda m_c - \frac{\lambda}{2r(1-r)^2}\label{g1a}\\
d_1 &=& g_1 + g_2 r + \frac{\lambda(r-2)}{2(r-1)}\label{d1a}
\end{eqnarray}
\end{flushleft}

The analytical solution of the difference equation is compared with the numerical solution obtained from FDM and FEM in fig. \ref{veraip1:a} \& \ref{veraip1:b}.





\begin{IEEEbiography}[{\includegraphics[width=1in,height=1.25in,clip,keepaspectratio]{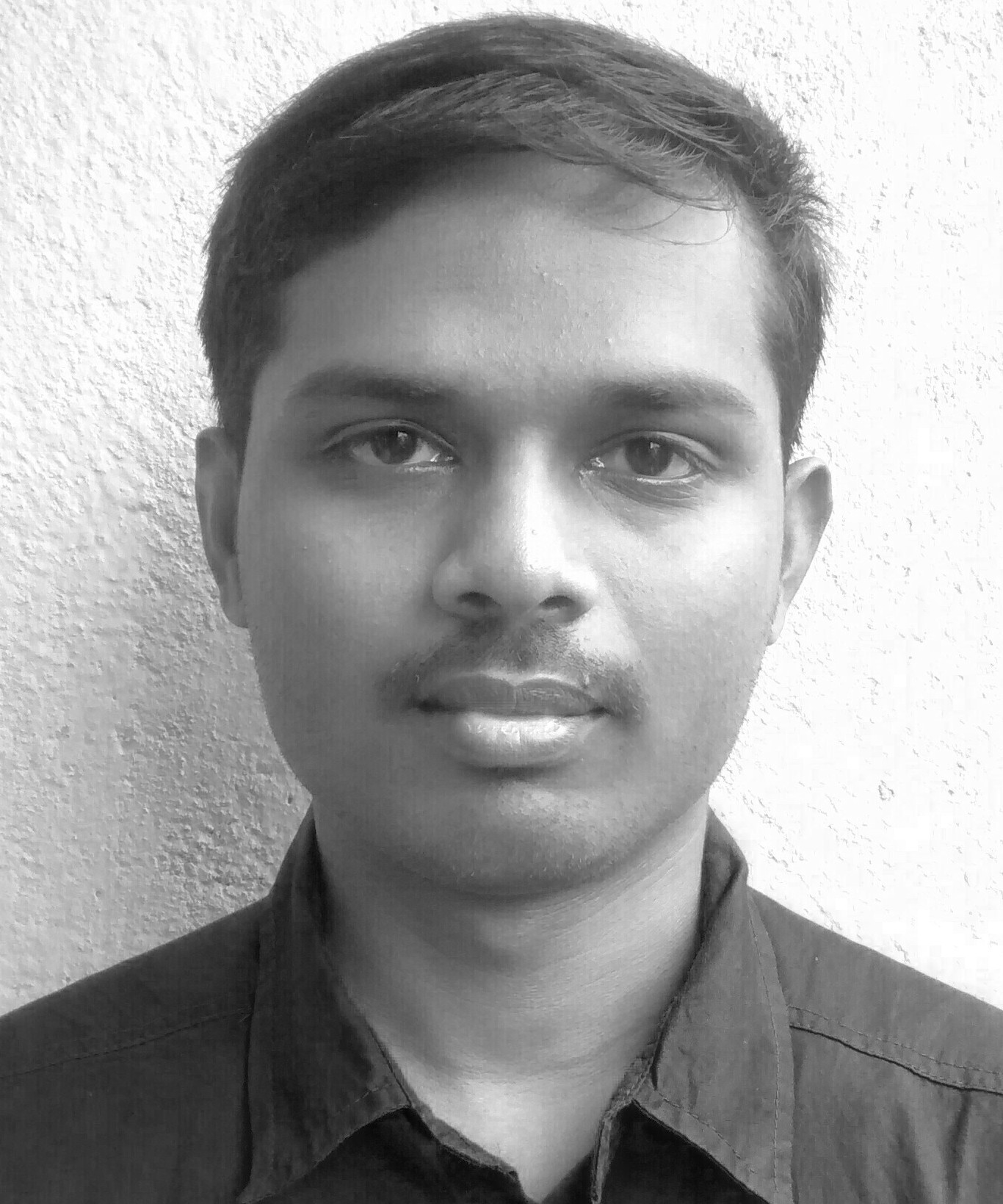}}]{Sethupathy S.}
received the B.E. degree in electrical and electronics engineering from the Anna University, Chennai, India, in 2009. He received the M.E. degree in electrical engineering from the Indian Institute of Science, Bangalore, India in 2011, where he is currently pursuing his Ph.D. degree. 

His research interests include Electromagnetism, Magnetohydrodynamics and Computational Methods.

\end{IEEEbiography}

\begin{IEEEbiography}[{\includegraphics[width=1in,height=1.25in,clip,keepaspectratio]{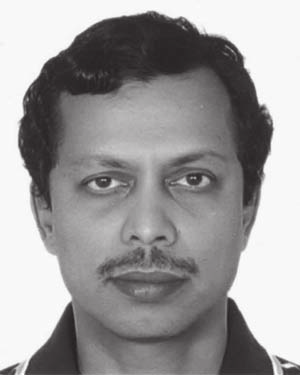}}]{Udaya Kumar}
 was born in Udupi District, Karnataka, India, in 1966. He received the bachelor’s degree in electrical engineering from Bangalore University, Bangalore, India, in 1989 and the M.E. and Ph.D. degrees in high-voltage engineering from the Indian Institute of Science, Bangalore, in 1991 and 1998, respectively. 

He was a Senior Analyst with the Electromagnetic Group of Electromagnetic Research Consultants, Bangalore. Since 1998, he has been with the Indian Institute of Science, where he is currently a Professor in the Department of Electrical Engineering. His research interests include lightning, electromagnetism, and high-frequency response of windings.

He is a member of CIGRE working groups WG C4.26 on “Evaluation of Lightning Shielding Analysis Methods for EHV and UHV DC and AC Transmission Lines” and WG C4.37 on “Electromagnetic Computation Methods for Lightning Surge Studies with Emphasis on the FDTD Method.”

\end{IEEEbiography}

\end{document}